\documentclass{article}

\usepackage{arxiv}

\usepackage[utf8]{inputenc} % allow utf-8 input
\usepackage[T1]{fontenc}    % use 8-bit T1 fonts
\usepackage{hyperref}       % hyperlinks
\usepackage{url}            % simple URL typesetting
\usepackage{booktabs}       % professional-quality tables
\usepackage{amsfonts}       % blackboard math symbols
\usepackage{nicefrac}       % compact symbols for 1/2, etc.
\usepackage{microtype}      % microtypography
\usepackage{lipsum}
\usepackage{amsmath}
\usepackage{graphicx}
\usepackage{enumerate}
\usepackage{lineno}
\usepackage{setspace}
\usepackage{amssymb}
\usepackage{verbatimbox}
\usepackage[english]{babel}
\usepackage{amsthm}

\title{a non-terminating game of beggar-my-neighbor}

\author{
  Brayden Casella \\
  Independent Researcher \\ Laconia, New Hampshire
\AND
    Philip M. Anderson \\
    Independent Researcher \\ Sedgefield, United Kingdom
\And
    Michael Kleber \\
    Google
\And
  Richard P. Mann$^*$ \\
  School of Mathematics \\ University of Leeds\\
   $^*$corresponding author: R.P.Mann@leeds.ac.uk
  \And
    Reed Nessler \\
    Department of Physics and Astronomy \\ Texas A\&M University
    \And 
    William Rucklidge\\
    Google
    \And
    Samuel G. Williams\\
    Independent Researcher \\ Denver, Colorado
    \And
    Nicolas Wu \\
    Department of Computer Science \\Imperial College London
}

\begin{document}
\maketitle

\begin{abstract}
We demonstrate the existence of a non-terminating game of Beggar-My-Neighbor, discovered by lead author Brayden Casella. We detail the method for constructing this game and identify a cyclical structure of 62 tricks that is reached by 30 distinct starting hands. We further present a short history of the search for this solution since the problem was posed, and a record of previously found longest terminating games. The existence of this non-terminating game provides a solution to a long-standing question which John H. Conway called an `anti-Hilbert problem.'
\end{abstract}

\section*{Introduction}
Beggar-My-Neighbor (\textsc{bmn})\footnote{
Additional names for this game include Beggar \textit{Your} Neighbor, Strip Jack Naked, Drive the Old Woman to Bed, and many others.}
is a simple deterministic card game for two players. The rules, as described in the 1911 edition of \textit{Encyclopaedia Britannica} are as follows:

\emph{`An ordinary pack is divided equally between two players, and the cards are held with the backs upwards. The first player lays down his top card face up, and the opponent plays his top card on it, and this goes on alternately as long as no court-card appears; but if either player turns up a court-card, his opponent has to play four ordinary cards to an ace, three to a king, two to a queen, one to a knave \footnote{In this paper we use the more typical terminology for the court (or face) cards: Ace (A); King (K); Queen (Q) and Jack (J).}, and when he has done so the other player takes all the cards on the table and places them under his pack;\footnote{The player who takes a trick is also the one who lays down the first face-up card to begin the subsequent trick.} if, however, in the course of this playing to a court-card, another court-card turns up, the adversary has in turn to play to this, and as long as neither has played a full number of ordinary cards to any court-card the trick continues. The player who gets all the cards into his hand is the winner.'}

Anyone playing \textsc{bmn} repeatedly will notice that some games are much longer than others. It is a natural question to ask whether there is a game which can continue indefinitely, and has surely occurred independently many times, but notably by John H. Conway during his time as a student at Cambridge \cite{paulhus1999}. Since there are only finitely many ways to arrange a deck of cards, a game continues forever if and only if at some point the game returns to exactly a situation it has been in before, i.e. game play enters a loop.
Conway further wrote,  ``The question of whether it always terminates is one
of my `anti-Hilbert' problems.  [Hilbert's were some problems for
mathematicians to work on during this century.  Mine are problems 
we should NOT be working on during the next!]''\footnote{We know of no definitive list of Conway's anti-Hilbert problems, but his Thrackle Conjecture is a second example; see \cite{stackexchange-antihilbert}}
\cite{mathfun1998, paulhus1999, nowakowski2002, gowers2008}.

Despite this inauspicious framing, a number of individuals (including all the authors of this paper) have since searched for a non-terminating game, and in the process have discovered games of increasing length. In this paper we give a short historical perspective on the search and demonstrate the existence of a non-terminating game.

\section*{History of the search}
The computer-assisted hunt for a non-terminating game of \textsc{bmn} first appears in print in \cite{beasley1989}.  Beasley reports on the statistical results of ten thousand randomly-dealt games, and notably observes that the game's length seems to have a half-life: ``if a game is down to two players but has not yet finished, the probability is about even that it will still not have finished within a further twenty tricks.''  He further notes that if the game were instead a \textit{random} map on possible game states with a 20-trick half-life, then finding a loop in the game state would be like looking for a needle in a haystack: a crude probabilistic argument says that there is over a 90\% chance of a loop existing, but the chance that a randomly-chosen starting position leads into a loop is around $1/10^{22}$.

In 1991, inspired by (or perhaps ``not deterred by'') this needle-in-a-haystack assessment, Rutgers University student Chris Long posted to the usenet newsgroup rec.puzzles, kicking off what he called \emph{The Great Usenet BYN Contest}.\cite{chrislong1991-start}  Long inaugurated the contest with the announcement that after a computer search of around 50 million random games, he had not found any that ran forever, but had found a game which ended after 538 tricks, and he challenged rec.puzzles contributors to find an infinite game or, failing that, to beat his record.  Long also reported that ``Jim Propp has informed me that John Conway has offered \$100 to anyone who can exhibit an initial BYN hand that has an infinite cycle.''  The contest ran for six weeks, and was won by Mark S. Manasse of Digital Equipment Corporation, who found a 713-trick game \cite{chrislong1991-end}, probably having examined on the order of 100 billion random games.

This unsolved problem gained greater prominence in 1999 with the publication by Marc Paulhus of an investigation into the structure of \textsc{bmn}, including the construction of non-terminating games in reduced sets of cards \cite{paulhus1999}, and the presentation of early records for the longest known terminating games of \textsc{bmn} with a full deck. Published by \emph{The American Mathematical Monthly} as part of its Unsolved Problems section, this brought the problem and its origins into wider visibility in the math community. 

An exhaustive test of all possible deals is computationally implausible: there are $\binom{52}{4}\binom{48}{4}\binom{44}{4}\binom{40}{4} \simeq 6.54 \times 10^{20}$ possible distinct starting configurations (ignoring the details of a deck of cards irrelevant to the game). As will be seen below, even very substantial computing resources can only hope to explore a small part of this space. Therefore, three approaches have suggested themselves: (i) statistical optimisation; (ii) brute force random search; and (iii) direct construction of a non-terminating game.

\begin{figure}
    \centering
    \includegraphics[width = 17cm]{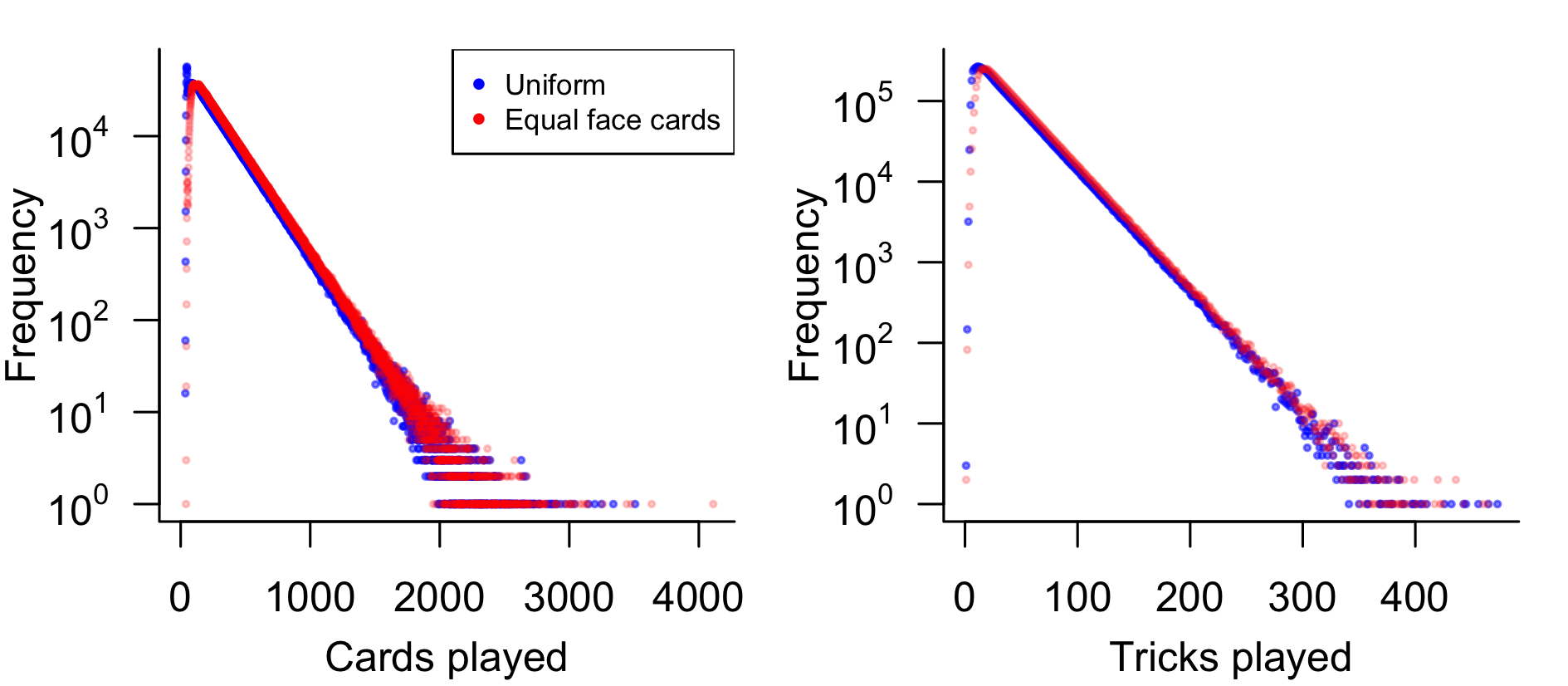}
    \caption{Frequency distribution of game lengths from $10^7$ randomly-dealt games, either with uniformly random deals (blue points) or restricted to deals in which both players have an equal share of the face cards (red points.) Apart from very short games, the frequencies in all cases closely resembles an exponential distribution (note the log-scale on the y-axis). Although dealing an equal share of face cards marginally increases expected game length, it does not change the decay rate of the exponential tail of the distribution.}
    \label{fig:stats}
\end{figure}

It is straightforward to computationally simulate many randomly-dealt games of \textsc{bmn} and to record the number of cards and tricks played during these games. Figure \ref{fig:stats} shows the frequency of different game lengths in a sample of $10^7$ randomly-dealt games. This closely resembles an exponential distribution, except in the case of very short games. The apparently-exponential distribution of game lengths suggests that although the progress of the game is deterministic, from a coarse-grained perspective (one in which we do not see the details of the packs) any ongoing game has a constant probability of `collapsing' as it progresses, and thus the fact that a game has endured for a certain amount time has no predictive power over its future length. The tail of this distribution closely matches Beasley's earlier observation of a half life of $\sim$20 tricks.

A statistical approach seeks to identify coarse-grained characteristics of starting packs that are associated with greater game length. An intuitively-sensible strategy is to ensure that each player is dealt the same number of each face card, such that the game is more evenly-balanced. As shown in Figure \ref{fig:stats} this strategy does have the effect of increasing the mean length of randomly-generated games. However, this is mostly achieved by selectively avoiding the very shortest games, and the exponential tail of the distribution appears to decay at the same rate as when games are dealt uniformly at random from all possible starting packs. As such this strategy can only very marginally increase the expected size of the longest game found over $n$ randomly-generated examples, and in our experience is an ineffective method for discovering longer games in the far tails of the distribution. Similar statistical approaches such as identifying specific individual card placements associated with higher mean game lengths also proved unsuccessful at systematically generating games of greater extreme length.

A large part of the search was subsequently carried out by brute force: simulating many randomly-dealt games in the hope of finding one that did not terminate, or as a lesser goal identifying one longer than any yet discovered. A significant innovation was to take found long games and play them `backwards', that is to identify earlier card configurations that result in the same configuration as those starting games do, such that this could represent an alternative starting position that produces a longer game. Although playing a game forwards is an entirely deterministic process, playing backwards involves some nondeterminism; families of different starting hands are able to converge on a particular state. These cases arise when looking through the history of a player's hand. When there is a matching amount of number cards after a face card, this could be because they won a turn, but it might also be due to some residual play. An example of games which lead to the same state can be seen in \hyperref[Appendix B]{Appendix B}. Hence, one must explore and keep track of an expanding tree of potential former states until none of these can be played backwards any further.

Nonetheless, despite this innovation, the search remained primarily constrained by the available computing power and time. To illustrate the scale of this search, Paulhus played $3.2 \times 10^9$ games prior publishing his initial findings \cite{paulhus1999}; Mann \& Wu (\hyperref[Appendix A]{Appendix A}) played c. $10^{13}$ games to find a record in 2007, running as a background process on between 5 and 10 desktop workstations in the Department of Engineering Science at the University of Oxford over a period of several weeks; Rucklidge (\hyperref[Appendix A]{Appendix A}) played $3.3 \times 10^{15}$ games in 2014 using 1000--8000 idle cores in Google's fleet, at a rate of about $10^9$ games per hour per core; Nessler played $3.5\times 10^{15}$ games over several years on a disused workstation capable of $7\times 10^{10}$ games per hour. 

\begin{figure}
    \centering
    \includegraphics[width=17cm]{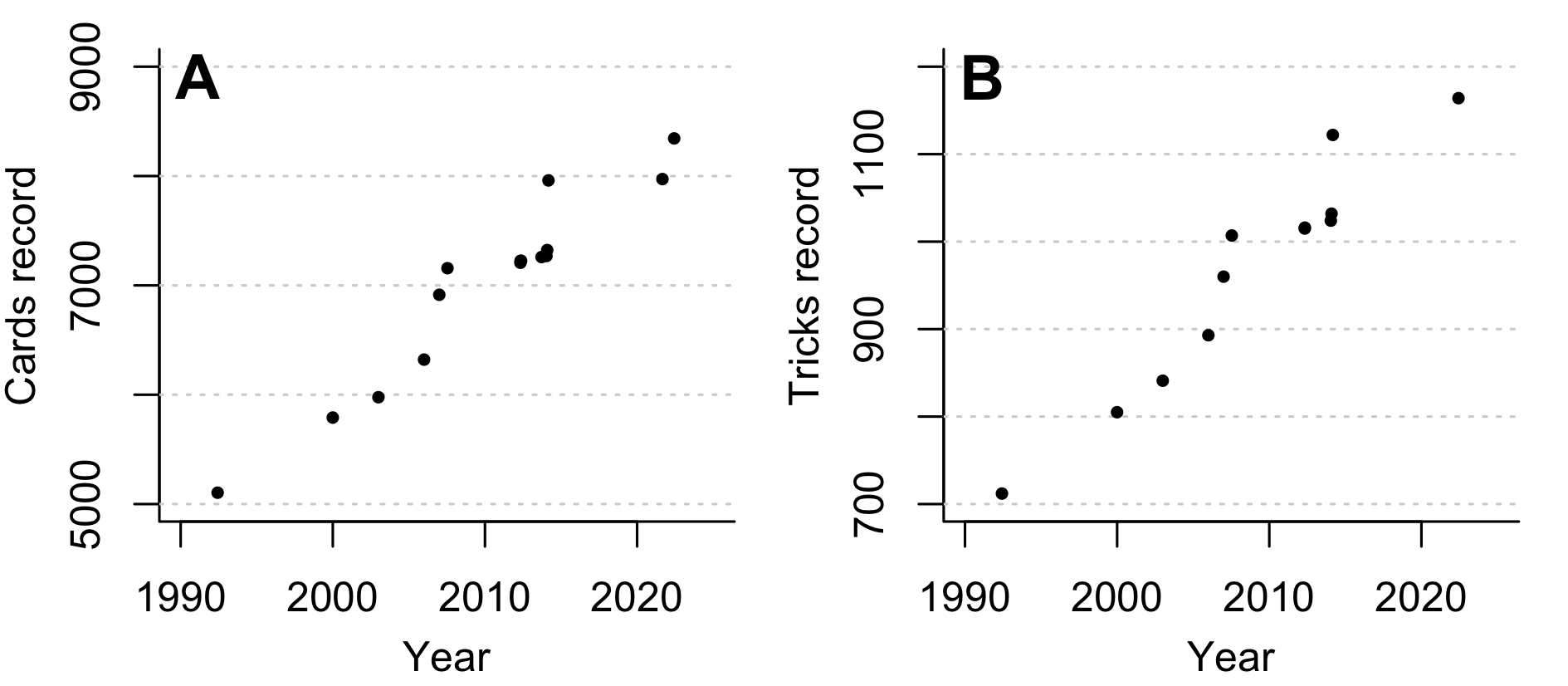}
    \caption{Historical records since 1992 for the longest known terminating game of \textsc{bmn} in terms of cards played (A) and tricks (B). A full list of historical records known to the authors is presented in \hyperref[Appendix A]{Appendix A}.
    }
    \label{fig:records}
\end{figure}

Figure \ref{fig:records} shows the historical progress of the records for known longest games, in terms of either cards played (A) or tricks (B). We stress that these are records known to us, which have been recorded via personal correspondence. Despite some stops and starts, progress in this record was broadly linear in time over the long term. This is straightforwardly explicable in terms of a simple model of random search. First, let us assume that games lengths are exponentially distributed with mean $\mu$, as indicated in Figure \ref{fig:stats} (this can apply to either cards or tricks). Assuming that a large number ($n$) of random games have been played, the expected length of the longest game ($r$) is given via standard order statistics as:
\begin{equation}
    \langle r \rangle = \mu \sum_{i=1}^n \frac{1}{n-i+1} = \mu H_n \simeq \mu (\gamma + \ln(n)), 
\end{equation}
where $H_n$ is the $n$th harmonic number and $\gamma$ is the Euler-Mascheroni constant. If we further assume that computing power increases exponentially over time and making the highly-simplifying assumption that attention to solving the problem remains constant (i.e. a constant fraction of available compute is dedicated to the task), then:
\begin{equation}
    n = Ae^{kt},
\end{equation}
for some constants $A$ and $k$. Substituting this value of $n$ into the earlier expression gives:
\begin{equation}
    \langle r \rangle = \mu (\gamma + \ln A + kt),
\end{equation}
showing that the longest game is expected to grow linearly with time $t$. 

\section*{Constructing a non-terminating game}
Playing with a standard deck of cards can be time consuming and difficult to gain understanding of the game. When the goal is to find a non-terminating game, it is also unclear how looking at terminating games will give insight into non-terminating games. Therefore, a good starting point is to find non-terminating games with the fewest number of cards.

Many smaller non-terminating decks with different distributions are easy to find. Consider two starting hands of a 6-card game: \verb|J--| and \verb|-J-| where the first hand begins play (hereafter notated as \verb|J--|/\verb|-J-|). The starting trick will transfer \verb|J-| to the first hand, then infinitely add \verb|-J-| to alternating hands. More complex solutions can be found through exhaustive search of smaller decks relatively quickly. Though solutions have been found using random search for decks up to 42 cards with various numbers of face cards, solve time grows rapidly, and brute force begins to break down.

\subsection*{Expansion, mutation, and the first non-terminating but unbalanced game}

Due to the repeating sequence of tricks in the \verb|J--|/\verb|-J-| game, it is possible to expand this simple deck to a slightly larger non-terminating game. Because the added tricks repeat, attaching those tricks to the end of the initial hands preserves the repeated sequence, and keeps the game non-terminating. In this case, it is done by doubling each player’s hand which makes the new game \verb|J--J--|/\verb|-J--J-|. More complicated games will have a longer repeating sequence of added tricks and a different way to expand into larger games. More generally, playing any non-terminating game can be thought of as a list of states, with each
state consisting of both player's hands:

\[ 
I_1 , I_2 \ldots , I_m , L_1, L_2, \ldots , L_n , L_1 , L_2 \ldots L_n , \ldots
\]

Here, each successive state corresponds to one
after a trick is added to a hand.
Each state $I_i$ where $i \in {1 \ldots m}$ corresponds to
an initial state before the loop,
and each state $L_j$ where $j \in {1 \ldots n}$ corresponds to
a state that is in the loop, which has a periodicity of $n$ states.
Each state eventually leads to a non-terminating game,
but only each $L_j$ will be repeated infinitely, where
the repeated states are:
\[ L_1, L_2 , \ldots , L_n \]

For each repeated state, player 1's hands can be concatenated to create a much larger starting hand, and the same done with player 2's hands. This new game will also be non-terminating. An additional step can be taken to expand further, where this expanded game can be attached any amount of times to make another non-terminating game for each addition. The decks created from this method are often not balanced (games where the length of the starting hands of each player are equal), but cutting the hands at equal length does often yield a non-terminating game.

This way, finding small non-terminating games with at least one of each face card might be able to be expanded into a standard deck which is non-terminating. To generate smaller games, an exhaustive search of balanced games up to 26 cards and various amounts of each face card can be generated and non-terminating games saved. A half-sized deck with half the face cards might be simply expanded into a standard deck, but this is not the case as there is no non-terminating game with this distribution \cite{paulhus1999}.

Forcing the games to be balanced is a severe limitation. To move into the unbalanced space, the thousands of saved non-terminating games can be expanded using the method described above. Each possible combination within those expansions of unbalanced games also tested for termination. Filtering the results for non-terminating decks which contain 4 of each face card eventually yields this deck:\\

\begin{verbatim}
 --Q-K-J-----------------KJ--Q----A-Q-A-A---J--Q-K-A-K
 J-
\end{verbatim}

There are 39 instead of 36 number cards but the three number cards in the \verb|A---J| sequence can be removed to make a non-terminating standard deck. This game does not become balanced at any state in the loop. Moreover, playing backwards to explore all possible earlier states that lead to the loop yields no balanced deck. A different deck is needed to construct a cycle that is accessible from an initial deal.

Small changes to non-terminating games can often be found which create new non-terminating games.
Three operations can be used to do this. For the first, at every point in the deck, each type of card (number card, jack, queen, king, and ace) is temporarily placed into the deck at that spot, then the game is played and tested for termination. If the game is non-terminating, that new deck is saved. The second operation is removal of cards, where each card in the deck is temporarily removed, then the deck tested and saved if non-terminating. The final operation is swapping one card for another. Similarly to the previous two operations, each card in the deck is temporarily changed to every other card and each change tested for termination. All points to split the deck between player 1 and player 2 are also tested for each operation. The speed of these operations is an issue so only a maximum of three additions and removals at a time were performed. Through the operations and manual changes, this non-terminating standard deck can be found:
\begin{verbatim}
 --------------------J--Q-Q-Q-Q-K-KJ--K-K-A-A-A-AJ-
 J-
\end{verbatim}

At no point in the state loop of this game do the hands become balanced, nor do any backtracked predecessors of the loop.

It is very difficult to create further non-terminating standard games using the operations and manual changes from the two games found. A tree search through the operations also failed to yield any interesting non-terminating games.

\subsection*{Deck structure and Casella's Ansatz}

Analyzing the two non-terminating standard games reveals why no state in a loop becomes close to balanced and brings about the key to eventually finding a balanced game.
Both non-terminating standard games continually cycle through adding sequences terminated by the jack cards to the end of player 1's hand before moving onto the next sequence. For example, the second non-terminating standard game shown earlier is this:
\begin{verbatim}
 --------------------J--Q-Q-Q-Q-K-KJ--K-K-A-A-A-AJ-
 J-
\end{verbatim}

Selecting specific states (in order, each 4 tricks apart) from the full trace of this game, we observe:
\begin{verbatim}
 --------------------J--Q-Q-Q-Q-K-KJ--K-K-A-A-A-AJ- / J-

 --Q-Q-Q-Q-K-KJ--K-K-A-A-A-AJ--------------------J- / J-

 --K-K-A-A-A-AJ--------------------J--Q-Q-K-Q-K-QJ- / J-

 --------------------J--Q-Q-K-Q-K-QJ--A-K-A-K-A-AJ- / J-
\end{verbatim}

Notice how player 1's hand mostly looks as though the deck is only being `slid' over until it reaches a jack. Sliding this way three times will result in an almost identical deck to the original, where all that is different is the order of a few face cards. Continuing play will return the face cards to the original order. The full trace of the solution can be found in \hyperref[Appendix D]{Appendix D}, where this can also be seen. Some different sequences will not be added to the end of player 1's hand as it was originally, but instead have their own loop and will reach the original when played through enough. Importantly, player 2's hand will always return to \verb|J-| before moving on to play through the next sequence in player 1's hand. These jack-terminated sequences, or `pieces' as they will be referred to, are therefore independent from one another and can be interchanged in any order to create a non-terminating deck. Due to this, player 2's hand can be seen only as a place to temporarily hold part of the current piece being played.

Inspired by this observation, we consider this promising-seeming shape for more non-terminating decks:\\

\newtheorem*{CA}{Casella's Ansatz}
\begin{CA}To search for an infinite game of \textsc{bmn}:
\begin{enumerate}
    \item 
    Look for games which pass through the state \\
    
    $\Delta_1$ \verb|J| $\Delta_2$ \verb|J| $\Delta_3$ \verb|J-|\\
    \verb|J-|\\
    
    where deck pieces $\Delta_1$, $\Delta_2$, and $\Delta_3$ behave independently, and game play proceeds by cycling these three deck pieces after $t_1$, $t_2$, and $t_3$ tricks respectively.  Some small internal permutations of face cards within each $\Delta_i$ may occur without disrupting the overall pattern.
    \item 
    To find a balanced starting position with 26 cards in player 2's hand, look for:
    \begin{enumerate}
        \item A deck piece $\Delta_i$ which can be formed by forward play from a state with more cards in player 2's hand (which we may search for by playing pieces backwards), and
        \item A deck piece $\Delta_j$ which moves a lot of cards over to player 2's hand in the middle of its forward play (so that the balanced initial state will represent the middle of playing through $\Delta_j$).\\
    \end{enumerate}
\end{enumerate}
\end{CA}

The precise non-terminating games found by the search up to this point are no longer relevant except as inspiration for this new method.  Finding the correct combination of pieces is the primary focus and what leads to the balanced non-terminating standard deck.

With the two decks currently available, the pieces are: \verb|--Q-Q-Q-Q-K-K|, \verb|-----------------K|, \verb|--Q----A-Q-A-A---|, and \verb|--------------------|. The first deck found appears to add an additional two pieces but that is due to the state being part way through dismantling and reconstructing its \verb|--Q-K-A-KQ-K-| piece. In these pieces, the most important aspect is the amount of face cards. Since number cards can often be added and removed from the same piece, and face cards can often be swapped for any other face card, the amount of strictly unique pieces can become overwhelming.
The \verb|--Q-Q-Q-Q-K-K| piece in particular can have any face card swapped for any other, while also able to add and remove various number cards. So

\verb|--Q-Q-Q-Q-K-K|,\\
\verb|--Q-Q-Q-Q-K-A|, \\
\verb|--Q-Q-Q-Q-K-A-|, and\\
\verb|--Q-Q-Q-Q-KK|

are all different pieces, but have the same amount of face cards with only minor mutations. The number of different pieces generated through face card swapping is $3^6$ from the \verb|--Q-Q-Q-Q-K-K| piece. With a number card added there are another $3^6$ unique pieces added with face card swapping. Since there are many ways to add and remove number cards from this piece, there are thousands of pieces which are only slightly different. Saving all of these for manual construction would not be useful. Instead, A single one of these pieces can be saved and the operations used to find the variations. Due to this, pieces which only differ through the operations (or repeated steps through the operations) are considered the same piece. For example, \verb|--Q-Q-Q-Q-K-K| and \verb|--Q-K-A-KQ-K-| are seen as the same piece.

The realization of the pieces also results in a natural template, or ‘skeleton’, deck for checking if a sequence is a piece. Since all pieces have a structure which creates independence between the pieces in a game, any sequence which doesn't have this structure will likely not create a non-terminating game (and therefore not be a piece). The structure of the template used to test for pieces is:\\

\verb| --J| FILTER \verb|J| TEST \verb|J-|\\
\verb| J-|

If no filter is used in the template, there will be sequences which lead to non-terminating games in the template but do not act independently of other pieces when used in a deck together (and therefore don't create a non-terminating game when used with other pieces). 
These false positives can be easily rejected by adding a properly selected found piece as the filter. The piece \verb|--K---A----AA| is a suitable filter and was selected due to its beneficial structure for backwards play (will be expanded upon later). The template is now:\\

\verb| --J--K---A----AAJ| TEST \verb|J-|\\
\verb| J-|

Any sequence can replace TEST in this template, just as \verb|--K---A----AA| replaced FILTER. If the game is non-terminating with this sequence, then the sequence is a piece and can be saved. The template arises from the structure of the pieces, and any amount of pieces (each terminated by a jack) can be added to the start of player 1's hand to create a non-terminating game.

This template based on jacks can be substituted for one based on the other face cards. However, this requires the inclusion of jacks in the inserted pieces and when this is done, one finds that many number cards are required in those pieces, which forces any constructed deck using these pieces to have an excess. As such, the jack template is used here.

Testing all combinations of a standard deck is an exhaustive search of the entire space, which is too large for an exhaustive search, but using templates this search space is vastly reduced.
The search space is not only reduced by the cards inherent to the template (four jacks and two number cards), but also due to a piece reducing the number of cards needed for the other pieces used in a deck. Therefore, all combinations of increasing size are tested with the template. To quickly find small pieces, number cards are considered 0, queens 1, kings 2, and aces 3. Finding all possible combinations is the same as counting in base 4 and padding with a sufficient amount of zeros (a total of 40 digits were used). While testing these, the leftmost digit can be repeatedly taken off the sequence and the reduced sequence tested again until a face card is reached, then move onto the next 40 card sequence. While this works exactly as desired early on, once more than 4 of a single face card or more than 12 total face cards start being generated, those sequences being tested are invalid for a standard deck of cards. As more digits are being populated by face cards, the time it takes to get to the next digit will also take exponentially more time. While that stays true, multiset permutations can be used to only sequences with a specific amount of each card, speeding up the testing significantly. Another speedup here can be made by assuming in most pieces with enough face cards, there will be face cards which can be replaced with any other face card. For example, instead of testing sequences that include 2 aces, 4 kings, and 4 queens, testing sequences with 6 kings and 4 queens works very well in reducing the time it takes to find larger pieces. Then when a piece is found, checking replacement of the extra kings into aces can be done.

Ten different non-terminating standard decks, each part of different state loops, can be quickly constructed with the new pieces and an eventual 15 being created before the final game.
Attempting backwards play at each state in the loop for all of these yields unbalanced games.
All pieces and non-terminating standard games constructed from them can be found in \hyperref[Appendix C]{Appendix C}.

\subsection*{Playing backwards to find balance}

Since the main focus is creating a balanced game, pieces need to be used which add as many cards to player 2's hand as possible. This can only be done through normal and backwards play of each piece. Beneficial backwards play for each piece mainly relies on subsequences within the piece. Looking back at the \verb|--K---A----AA| piece used as the filter in the template, completing this with jack termination creates the non-terminating game:
\begin{verbatim}
 --K---A----AAJ-
 J-
\end{verbatim}

To play backwards, what do the starting hands need to be to add an equivalent subsequence of cards to the end of player 1's hand? Looking from right to left of player 1's hand, \verb|-| is not a valid trick (because play would continue normally rather than the card added to a player's hand), so the first valid subsequence is \verb|J-|. To add \verb|J-| as a trick through normal play, the game would be \verb|J|/\verb|-|. These cards can be removed from the end of player 1's hand and the game which plays \verb|J-| added to the start of the original hands. This creates:
\begin{verbatim}
 J--K---A----AA
 -J-
\end{verbatim}

Playing this game forward one trick will create the original deck. If the same backwards play is attempted for player 1's hand now, the effort is immediately stopped since no valid trick ends with \verb|A| in normal play. Therefore, enough number cards are needed to be left for the next face card before the subsequence being removed. Going back to the original deck, the subsequence which does this is \verb|AAJ-|. Playing this backwards creates:
\begin{verbatim}
 AJ--K---A----
 A-J-
\end{verbatim}

Notice how four cards were removed from the end of player 1's hand, then two cards added to the start of each hand. The beneficial structure of the selected piece shows itself once backwards play is continued from here. The subsequence \verb|A----| is selected to be removed from player 1's hand and played backwards, creating:
\begin{verbatim}
 AAJ--K---
 ----A-J-
\end{verbatim}

Five cards are removed while only adding one to the start of player 1's hand, and four to player 2's hand. When the goal is to make player 2's hand larger, this is an extremely efficient subsequence. If the \verb|K---| subsequence is selected for backwards play next, another four cards are removed from player 1's hand while only adding one card to the start of their hand. Not every piece moves so many cards into player 2's hand, and some pieces have subsequences which swap whose turn it is to play after playing backwards. This is disastrous if it is not the first subsequence in the piece and results in the entire piece being unusable to create a balanced deck (for all pieces found). As an example, this is a promising looking piece where the necessary \verb|-KK--Q-K--Q--| subsequence in the piece swaps whose turn it is after backwards play:\\

\verb| --------K----KK--Q-K--Q--QQ|\\

Now, all that’s needed is a piece which allows for a balanced game. Eventually, this piece is discovered:\\ 

\verb| --------K---------Q-Q-K---Q-KKQ|\\

Testing the piece shows that two of the face cards can be replaced with any other face card. Manually selecting the pieces \verb|--K---A----AA| and \verb|--|, then adding all three together with jack termination creates this deck:\\

\begin{verbatim}
 --K---A----AAJ--J--------K---------Q-Q-K---Q-KAQJ-
 J-
\end{verbatim}

Playing forward naturally through the first piece, then backwards for the final piece yields the balanced non-terminating standard game:
\begin{verbatim}
 ---K---Q-KQAJ-----AAJ--J--
 ----------Q----KQ-J-----KA
\end{verbatim}

\section*{The first non-terminating game and its family}
The non-terminating game constructed by the method above is given by the starting hands below, with player 1 being the starting player, and with the top of the hand being to the left.\\

Player 1: \verb|---K---Q-KQAJ-----AAJ--J--| \\
Player 2: \verb|----------Q----KQ-J-----KA| \\

This starting deal leads after 4 tricks to a cycle of 62 tricks of which the first repeating state is:\\

Player 1: \verb|---AAJ--J--------K---------Q-Q-K---Q-KAQJ-| \\
Player 2: \verb|-----KA-J-|\\

with player 2 to play next. The full progress of the cycle is given in \hyperref[Appendix D]{Appendix D}.

\begin{figure}[h!]
    \centering
    \includegraphics[width=15cm]{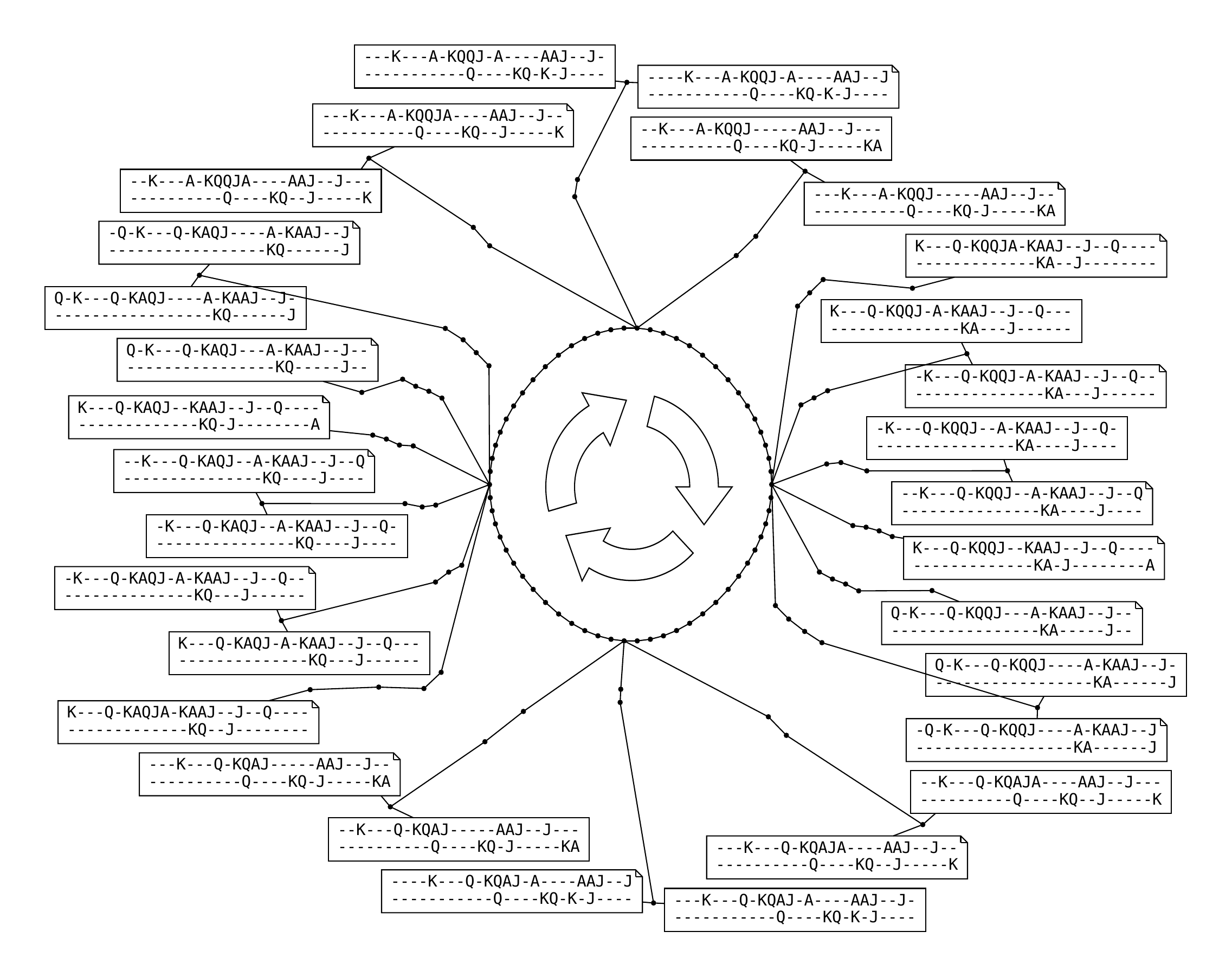}
    \caption{Illustration of a family of non-terminating games, all converging to the same 62-trick cycle. Starting packs shown in dogeared boxes are those with the topmost player starting; those in plain boxes start with the lower player.}
    \label{fig:family}
\end{figure}

By playing backwards from each state within the cycle, it is possible to identify a family of 30 possible starting deals, each of which leads to a point on this cycle and thus each of which represents a non-terminating game. These are illustrated in Figure \ref{fig:family}.

\section*{Conclusion}
We have demonstrated the existence of a non-terminating game of Beggar-My-Neighbor, providing a solution to the `anti-Hilbert problem' posed by John H. Conway. Having identified a cyclical game, we were able to further construct a family of related solutions that converge to the same cycle. Although the original open problem is now resolved, several questions remain open. First, how many possible cycles can be reached? Second, how many balanced non-terminating games exist? Third, what is the length of the longest terminating game? Before the construction of this example, in excess of $10^{15}$ games were simulated without generating a non-terminating game, so players can rest assured that the probability of encountering such a case `in the wild' are negligible.

% Bibliography
\bibliographystyle{ieeetr}

\clearpage

\section*{Appendix A: Historical record of longest known terminating games}
\label{Appendix A}
Note: these are records known to the authors and established by personal correspondence. Each record shows the starting packs, with cards on the left representing the top of each pack. The topmost player starts.

\begin{myverbbox}{\manasse}
 Q--J----K--K-J---Q---A---A
 ---K-K--JA-QA--J-----Q----
\end{myverbbox}
\begin{myverbbox}{\kleberA}
 ---JQ---K-A----A-J-K---QK-
 -J-----------AJQA----K---Q
\end{myverbbox}
\begin{myverbbox}{\kleberB}
 ----QJ-A-KK--------K--QQJ-
 ---J---K--A-----Q-AJ-A----
\end{myverbbox}
\begin{myverbbox}{\kleberC}
 -J--KA----A-Q--Q-A----KJ--
 A------QKJ--Q-------KJ----
\end{myverbbox}
\begin{myverbbox}{\collins}
 A-QK------Q----KA-----J---
 -JAK----A--Q----J---QJ--K-
\end{myverbbox}
\begin{myverbbox}{\mannwu}
 K-KK----K-A-----JAA--Q--J-
 ---Q---Q-J-----J------AQ--
\end{myverbbox}
\begin{myverbbox}{\nesslerA}
 ----Q------A--K--A-A--QJK-
 -Q--J--J---QK---K----JA---
\end{myverbbox}
\begin{myverbbox}{\nesslerB}
 -J-------Q------A--A--QKK-
 -A-Q--J--J---Q--AJ-K---K--
\end{myverbbox}
\begin{myverbbox}{\wu}
 --A-Q--J--J---Q--AJ-K---K-
 -J-------Q------A--A--QKK-
\end{myverbbox}
\begin{myverbbox}{\nesslerC}
 ---AK-Q--J----J--QKJ-Q----
 ------JK-----A--K--Q---AA-
\end{myverbbox}
\begin{myverbbox}{\rucklidgeA}
 A-AQ-----Q--K--AQ-------JJ
 -J-A-KKJ--K-----------Q---
\end{myverbbox}
\begin{myverbbox}{\anderson}
 -AJ--QK--K----Q--J-A-KKJ--
 ---------JQ----------A-AQ-
\end{myverbbox}
\begin{myverbbox}{\rucklidgeB}
 -J------Q------AAA-----QQ-
 K----JA-----------KQ-K-JJK
\end{myverbbox}
\begin{myverbbox}{\nesslerD}
 ----K---A--Q-A--JJA------J
 -----KK---------A-JK-Q-Q-Q
\end{myverbbox}
\begin{myverbbox}{\nesslerE}
 ---AJ--Q---------QAKQJJ-QK
 -----A----KJ-K--------A---
\end{myverbbox}

\renewcommand{\arraystretch}{2.5}
\begin{tabular}{ll}
Mark Manasse (1992) 713 tricks, 5104 cards \cite{chrislong1991-end}  & \parbox{6cm}{\manasse} \\
Michael Kleber (before 1999). 805 tricks, 5791 cards \cite{paulhus1999} & \parbox{6cm}{\kleberA}\\
Michael Kleber (c. 2002). 841 tricks, 5977 cards & \parbox{6cm}{\kleberB}\\
Michael Kleber (c. 2005). 893 tricks, 6321 cards & \parbox{6cm}{\kleberC}\\
Truman Collins (2006). 960 tricks, 6914 cards & \parbox{6cm}{\collins}\\
Richard Mann \& Nicolas Wu (16-07-2007). 1007 tricks, 7157 cards & \parbox{6cm}{\mannwu}\\
Reed Nessler (01-05-2012). 1015 tricks, 7207 cards & \parbox{6cm}{\nesslerA}\\
Reed Nessler (04-05-2012). 1016 tricks, 7224 cards & \parbox{6cm}{\nesslerB}\\
Nicolas Wu (17-05-2012). 1016 tricks, 7225 cards & \parbox{6cm}{\wu}\\
Reed Nessler (20-09-2013). 1014 tricks, 7259 cards & \parbox{6cm}{\nesslerC}\\
William Rucklidge (17-01-2014). 1024 tricks, 7269 cards & \parbox{6cm}{\rucklidgeA}\\
Philip Anderson (03-02-2014). 1032 tricks, 7323 cards & \parbox{6cm}{\anderson}\\
William Rucklidge (05-03-2014). 1122 tricks, 7960 cards & \parbox{6cm}{\rucklidgeB}\\
Reed Nessler (31-08-2021). 1106 tricks, 7972 cards & \parbox{6cm}{\nesslerD}\\
Reed Nessler (09-06-2022). 1164 tricks, 8344 cards & \parbox{6cm}{\nesslerE}\\
\end{tabular}

\section*{Appendix B: Example of games which lead to the same state after one trick}
\label{Appendix B}

Note: the number in parentheses is the player who starts play.\\

Initial game:\\
1. \verb|--------------------J--Q-Q-Q-Q-K-KJ--K-K-A-A-A-AJ-| 2. \verb|J-| (1)\\

Games after one trick of backwards play from the initial deck:\\
1. \verb|J--------------------J--Q-Q-Q-Q-K-KJ--K-K-A-A-A-A| 2. \verb|-J-| (1)

1. \verb|J--------------------J--Q-Q-Q-Q-K-KJ--K-K-A-A-A-| 2. \verb|A-J-| (2)

1. \verb|-J--------------------J--Q-Q-Q-Q-K-KJ--K-K-A-A-A| 2. \verb|A-J-| (1)

1. \verb|AJ--------------------J--Q-Q-Q-Q-K-KJ--K-K-A-A-| 2. \verb|-A-J-| (1)

1. \verb|-AJ--------------------J--Q-Q-Q-Q-K-KJ--K-K-A-| 2. \verb|A-A-J-| (2)

1. \verb|--AJ--------------------J--Q-Q-Q-Q-K-KJ--K-K-A| 2. \verb|A-A-J-| (1)

1. \verb|A-AJ--------------------J--Q-Q-Q-Q-K-KJ--K-K-| 2. \verb|-A-A-J-| (1)

1. \verb|A-AJ--------------------J--Q-Q-Q-Q-K-KJ--K-K| 2. \verb|--A-A-J-| (2)

1. \verb|-A-AJ--------------------J--Q-Q-Q-Q-K-KJ--K-| 2. \verb|K-A-A-J-| (2)

1. \verb|--A-AJ--------------------J--Q-Q-Q-Q-K-KJ--K| 2. \verb|K-A-A-J-| (1)

1. \verb|K-A-AJ--------------------J--Q-Q-Q-Q-K-KJ--| 2. \verb|-K-A-A-J-| (1)

1. \verb|K-A-AJ--------------------J--Q-Q-Q-Q-K-KJ-| 2. \verb|--K-A-A-J-| (2)

1. \verb|-K-A-AJ--------------------J--Q-Q-Q-Q-K-KJ| 2. \verb|--K-A-A-J-| (1)

\section*{Appendix C: Found pieces and non-terminating standard games from them}
\label{Appendix C}

Note: only pieces based on the definition provided are listed, there are thousands of strictly unique pieces.\\

Found pieces:
\begin{verbatim}
 --
 --Q-Q
 -----K
 --------
 --K---Q-KQQ
 --Q-Q-Q-Q-K-K
 --K---A----AA
 --------A-KAA
 --A-----A-A--A
 --A----Q----QQ
 --------A----AA
 --Q----Q----Q-Q
 --Q--Q-Q-QQ-Q-Q
 --------Q-Q-Q-Q
 --------Q----Q-Q
 --K-A--Q-K--Q--QQ
 --Q----A-Q-A-A---
 --------------K---
 --------Q----Q----
 -----------------Q-
 --------------------
 --------K----KQQ-KQQ
 --------K----KK--Q-K--Q--QQ
 --------K---------Q-Q-K---Q-KKQ
 --------------------A----------A---A
 --------------------Q----------Q----Q-Q
\end{verbatim}

Constructed non-terminating standard games:

\begin{verbatim}
 --------------------J--Q-Q-Q-Q-K-KJ--K-K-A-A-A-AJ- / J-
 -----------------KJ--Q-K-A-K-Q-KJ--Q----A-Q-A-A-J- / J-
 -----------------KJ--A-A-Q-A-A---Q-J--K---Q-KQK-J- / J-
 --------Q----Q-Q-J--------A----AA---J--K-KQ-KK-AJ- / J-
 --Q----Q----Q-Q-J--K---K-KAK--J--------A----AA--J- / J-
 -----------------Q-J--K---Q-KQQJ--K-K--A-A-A--A-J- / J-
 --------J--Q----A-Q-A-A-J--------K----KQ-Q-K-K-AJ- / J-
 -----A--J--Q----Q----Q-QJ--------K----KA--A-KK-AJ- / J-
 --------Q-------Q-J--Q-Q-K-K-AKJ--K---A----AA---J- / J-
 --------Q-------Q-J--Q-Q-K-K-K-KJ--A-----A-A--A-J- / J-
 --Q----Q----Q-Q-J-----KJ--------K----KA-A-K-A--AJ- / J-
 --A-Q-K-K-K-KJ--Q----A-Q-A-AJ--------------Q----J- / J-
 --------------Q----J--K---Q-KQKJ--A-A-K-A-A---Q-J- / J-
 --K---A----AA---J--J--------K----KA--Q-K--Q--QQ-J- / J-
 --A-A--J--A-A---J--------K---------Q-Q-K---Q-KKQJ- / J-
 --K---A----AAJ--J--------K---------Q-Q-K---Q-KAQJ- / J-
\end{verbatim}

\section*{Appendix D: A trace of the non-terminating cycle}
\label{Appendix D}

Note: for each configuration, the number in parentheses is the player who starts play, i.e. the winner of the previous trick.

\begin{enumerate}
\item
1. \verb|--K---A----AAJ--J--------K---------Q-Q-K---Q-KAQJ-| 2. \verb|J-| (1)
\item
1. \verb|K---A----AAJ--J--------K---------Q-Q-K---Q-KAQJ-| 2. \verb|--J-| (2)
\item
1. \verb|--A----AAJ--J--------K---------Q-Q-K---Q-KAQJ-| 2. \verb|--K-J-| (2)
\item
1. \verb|---AAJ--J--------K---------Q-Q-K---Q-KAQJ-| 2. \verb|-----KA-J-| (2)
\item
1. \verb|--J--------K---------Q-Q-K---Q-KAQJ--------A-KAAJ-| 2. \verb|J-| (1)
\item
1. \verb|J--------K---------Q-Q-K---Q-KAQJ--------A-KAAJ-| 2. \verb|--J-| (2)
\item
1. \verb|--------K---------Q-Q-K---Q-KAQJ--------A-KAAJ--J-| 2. \verb|J-| (1)
\item
1. \verb|------K---------Q-Q-K---Q-KAQJ--------A-KAAJ--J-| 2. \verb|--J-| (2)
\item
1. \verb|---K---------Q-Q-K---Q-KAQJ--------A-KAAJ--J-| 2. \verb|-----J-| (2)
\item
1. \verb|--------Q-Q-K---Q-KAQJ--------A-KAAJ--J-| 2. \verb|--------K-J-| (2)
\item
1. \verb|Q-K---Q-KAQJ--------A-KAAJ--J-| 2. \verb|-----------------KQ-J-| (2)
\item
1. \verb|-K---Q-KAQJ--------A-KAAJ--J--Q--| 2. \verb|--------------KQ-J-| (1)
\item
1. \verb|---Q-KAQJ--------A-KAAJ--J--Q----K---| 2. \verb|----------KQ-J-| (1)
\item
1. \verb|-KAQJ--------A-KAAJ--J--Q----K---------Q--| 2. \verb|-----KQ-J-| (1)
\item
1. \verb|AQJ--------A-KAAJ--J--Q----K---------Q----K---| 2. \verb|-KQ-J-| (1)
\item
1. \verb|--------A-KAAJ--J--Q----K---------Q----K---A-KQQJ-| 2. \verb|J-| (1)
\item
1. \verb|------A-KAAJ--J--Q----K---------Q----K---A-KQQJ-| 2. \verb|--J-| (2)
\item
1. \verb|---A-KAAJ--J--Q----K---------Q----K---A-KQQJ-| 2. \verb|-----J-| (2)
\item
1. \verb|KAAJ--J--Q----K---------Q----K---A-KQQJ-| 2. \verb|--------A-J-| (2)
\item
1. \verb|AAJ--J--Q----K---------Q----K---A-KQQJ--K---| 2. \verb|----A-J-| (1)
\item
1. \verb|AJ--J--Q----K---------Q----K---A-KQQJ--K---A----| 2. \verb|A-J-| (1)
\item
1. \verb|--J--Q----K---------Q----K---A-KQQJ--K---A----AAJ-| 2. \verb|J-| (1)
\item
1. \verb|J--Q----K---------Q----K---A-KQQJ--K---A----AAJ-| 2. \verb|--J-| (2)
\item
1. \verb|--Q----K---------Q----K---A-KQQJ--K---A----AAJ--J-| 2. \verb|J-| (1)
\item
1. \verb|Q----K---------Q----K---A-KQQJ--K---A----AAJ--J-| 2. \verb|--J-| (2)
\item
1. \verb|---K---------Q----K---A-KQQJ--K---A----AAJ--J-| 2. \verb|--Q-J-| (2)
\item
1. \verb|--------Q----K---A-KQQJ--K---A----AAJ--J-| 2. \verb|-----Q-K-J-| (2)
\item
1. \verb|-Q----K---A-KQQJ--K---A----AAJ--J-| 2. \verb|-K-J-----------Q--| (2)
\item
1. \verb|---K---A-KQQJ--K---A----AAJ--J-| 2. \verb|-----------Q----KQ-J-| (2)
\item
1. \verb|---A-KQQJ--K---A----AAJ--J--------K---| 2. \verb|----Q----KQ-J-| (1)
\item
1. \verb|QQJ--K---A----AAJ--J--------K---------A-Q-K---| 2. \verb|-KQ-J-| (1)
\item
1. \verb|--K---A----AAJ--J--------K---------A-Q-K---Q-KQQJ-| 2. \verb|J-| (1)
\item
1. \verb|K---A----AAJ--J--------K---------A-Q-K---Q-KQQJ-| 2. \verb|--J-| (2)
\item
1. \verb|--A----AAJ--J--------K---------A-Q-K---Q-KQQJ-| 2. \verb|--K-J-| (2)
\item
1. \verb|---AAJ--J--------K---------A-Q-K---Q-KQQJ-| 2. \verb|-----KA-J-| (2)
\item
1. \verb|--J--------K---------A-Q-K---Q-KQQJ--------A-KAAJ-| 2. \verb|J-| (1)
\item
1. \verb|J--------K---------A-Q-K---Q-KQQJ--------A-KAAJ-| 2. \verb|--J-| (2)
\item
1. \verb|--------K---------A-Q-K---Q-KQQJ--------A-KAAJ--J-| 2. \verb|J-| (1)
\item
1. \verb|------K---------A-Q-K---Q-KQQJ--------A-KAAJ--J-| 2. \verb|--J-| (2)
\item
1. \verb|---K---------A-Q-K---Q-KQQJ--------A-KAAJ--J-| 2. \verb|-----J-| (2)
\item
1. \verb|--------A-Q-K---Q-KQQJ--------A-KAAJ--J-| 2. \verb|--------K-J-| (2)
\item
1. \verb|Q-K---Q-KQQJ--------A-KAAJ--J-| 2. \verb|-----------------KA-J-| (2)
\item
1. \verb|-K---Q-KQQJ--------A-KAAJ--J--Q--| 2. \verb|--------------KA-J-| (1)
\item
1. \verb|---Q-KQQJ--------A-KAAJ--J--Q----K---| 2. \verb|----------KA-J-| (1)
\item
1. \verb|-KQQJ--------A-KAAJ--J--Q----K---------Q--| 2. \verb|-----KA-J-| (1)
\item
1. \verb|QQJ--------A-KAAJ--J--Q----K---------Q----K---| 2. \verb|-KA-J-| (1)
\item
1. \verb|--------A-KAAJ--J--Q----K---------Q----K---Q-KQAJ-| 2. \verb|J-| (1)
\item
1. \verb|------A-KAAJ--J--Q----K---------Q----K---Q-KQAJ-| 2. \verb|--J-| (2)
\item
1. \verb|---A-KAAJ--J--Q----K---------Q----K---Q-KQAJ-| 2. \verb|-----J-| (2)
\item
1. \verb|KAAJ--J--Q----K---------Q----K---Q-KQAJ-| 2. \verb|--------A-J-| (2)
\item
1. \verb|AAJ--J--Q----K---------Q----K---Q-KQAJ--K---| 2. \verb|----A-J-| (1)
\item
1. \verb|AJ--J--Q----K---------Q----K---Q-KQAJ--K---A----| 2. \verb|A-J-| (1)
\item
1. \verb|--J--Q----K---------Q----K---Q-KQAJ--K---A----AAJ-| 2. \verb|J-| (1)
\item
1. \verb|J--Q----K---------Q----K---Q-KQAJ--K---A----AAJ-| 2. \verb|--J-| (2)
\item
1. \verb|--Q----K---------Q----K---Q-KQAJ--K---A----AAJ--J-| 2. \verb|J-| (1)
\item
1. \verb|Q----K---------Q----K---Q-KQAJ--K---A----AAJ--J-| 2. \verb|--J-| (2)
\item
1. \verb|---K---------Q----K---Q-KQAJ--K---A----AAJ--J-| 2. \verb|--Q-J-| (2)
\item
1. \verb|--------Q----K---Q-KQAJ--K---A----AAJ--J-| 2. \verb|-----Q-K-J-| (2)
\item
1. \verb|-Q----K---Q-KQAJ--K---A----AAJ--J-| 2. \verb|-K-J-----------Q--| (2)
\item
1. \verb|---K---Q-KQAJ--K---A----AAJ--J-| 2. \verb|-----------Q----KQ-J-| (2)
\item
1. \verb|---Q-KQAJ--K---A----AAJ--J--------K---| 2. \verb|----Q----KQ-J-| (1)
\item
1. \verb|QAJ--K---A----AAJ--J--------K---------Q-Q-K---| 2. \verb|-KQ-J-| (1)
\end{enumerate}

\end{document}